\documentclass[12pt]{article}
\usepackage{amssymb,amsfonts,amsmath, psfrag,eepic,colordvi,graphicx,epsfig,ytableau}
\usepackage{amssymb,latexsym,graphics,array}
\usepackage{MnSymbol}
\usepackage[enableskew]{youngtab}
\usepackage{float,enumerate,enumitem}
\usepackage{tikz,ifthen}
\usepackage{pgfplots}
\usetikzlibrary{math}
\topskip  
\numberwithin{equation}{section}
\newtheorem{theorem}{Theorem}[section]

\newtheorem{definition}[theorem]{Definition}

\newtheorem{remark}[theorem]{Remark}
\newtheorem{lemma}[theorem]{Lemma}

\newtheorem{example}[theorem]{Example}

\newtheorem{problem}[theorem]{Problem}
\newtheorem{question}[theorem]{Question}

\begin{document}
	\parskip 6pt
	
	\pagenumbering{arabic}
	\def\sof{\hfill\rule{2mm}{2mm}}
	\def\ls{\leq}
	\def\gs{\geq}
	\def\SS{\mathcal S}
	\def\qq{{\bold q}}
	\def\MM{\mathcal M}
	\def\TT{\mathcal T}
	\def\EE{\mathcal E}
	\def\lsp{\mbox{lsp}}
	\def\rsp{\mbox{rsp}}
	\def\pf{\noindent {\it Proof.} }
	\def\mp{\mbox{pyramid}}
	\def\mb{\mbox{block}}
	\def\mc{\mbox{cross}}
	\def\qed{\hfill \rule{4pt}{7pt}}
	\def\block{\hfill \rule{5pt}{5pt}}
	\def\lr#1{\multicolumn{1}{|@{\hspace{.6ex}}c@{\hspace{.6ex}}|}{\raisebox{-.3ex}{$#1$}}}
	\def\red{\textcolor{red}}

	\begin{center}
{\Large \bf Difference ascent sequences and related combinatorial structures}
	\end{center}
	
	\begin{center}
		{\small Yongchun Zang$^a$,  Robin D.P. Zhou$^{a,b,*}$\footnote{$^*$Corresponding author.}  \footnote{{\em E-mail address:} dapao2012@163.com. }}
		
		$^a$College of Mathematics Physics and Information\\
		Shaoxing University\\
		Shaoxing 312000, P.R. China\\
		$^b$Institute of Artificial Intelligence\\
		Shaoxing University\\
		Shaoxing 312000, P.R. China
		
	\end{center}
	
\noindent {\bf Abstract.} 	
Ascent sequences were introduced by Bousquet-M\'elou, Claesson, Dukes and Kitaev, and are in bijection with unlabeled $(2+2)$-free posets, Fishburn matrices, permutations
avoiding a bivincular pattern of length $3$, and Stoimenow matchings.
Analogous results for weak ascent sequences have been obtained by
B\'enyi, Claesson and Dukes.
Recently, Dukes and Sagan introduced a more general class of sequences which are called $d$-ascent sequences.
They showed that some maps from the weak case can be extended to 
bijections for general $d$ while the extensions of others continue to 
be injective but not surjective.
The main objective of this paper is to restore these injections to bijections.
To be specific,  we introduce a class of permutations which we call  
difference $d$ permutations and a class of factorial posets which we call  difference $d$ posets, 
both of which are shown to be in bijection with $d$-ascent sequences.
Moreover, we also give a direct bijection between 
a class of matrices with a certain column restriction and 
Fishburn matrices.
Our results give answers to several questions posed by Dukes and Sagan.



	
	\section{Introduction}
	
	Let $x = x_1x_2\cdots x_n$ be a sequence of integers.
	An index $i$ ($1\leq i \leq n-1$) is said to be an {\em ascent} of $x$ if 
	$x_{i+1} > x_i$.
	Let $\mathrm{asc}(x)$ 
	denote the number of ascents of $x$.
	We call a sequence $x$ an {\em ascent sequence} if $x_1=0$ and 
	$0 \leq x_i \leq \mathrm{asc}(x_1x_2\cdots x_{i-1})+1$ for all $2 \leq i \leq n$.
	For example, one can check that $x = 01021324$ is an ascent sequence
	while $x = 0122431$ is not an ascent sequence.
	Let $\mathcal{A}_n$ denote the set of ascent sequences of 
	length $n$.
	For example, we have 
\[
\mathcal{A}_3=\{000,001,010,011,012\}.
\]

Ascent sequences were introduced by Bousquet-M\'elou et al. \cite{Bousquet}
to unify three other combinatorial structures: unlabeled $(2+2)$-free posets, permutations
avoiding a bivincular pattern of length $3$ and Stoimenow matchings \cite{Stoimenow}.
They have since evolved into a research hotspot, drawing considerable attention from scholars such as \cite{Chen2013,Claesson1,Claesson2,Conway,Dukes2010,Dukes2016,Dukes2019,Duncan,Fu,Jin,Lin,Mansour,Pudwell,Yan2014,Yan2018}.

B\'enyi et al. \cite{Benyi}  studied weak ascent sequences.
Given a sequence $x = x_1x_2\cdots x_n$, an index $i$ ($1\leq i \leq n-1$) is said 
to be a {\em weak ascent} of $x$ if $x_{i+1} \geq x_i$. 
Let $\mathrm{wasc}(x)$  denote the number of weak ascents of $x$.
The sequence $x$ is called a {\em weak ascent sequence} if  $x_1=0$ and 
$0 \leq x_i \leq \mathrm{wasc}(x_1x_2\cdots x_{i-1})+1$ for all $2 \leq i \leq n$.
Even though $x = 0122431$ is not an ascent sequence, it is a weak ascent sequence.
In the spirit of \cite{Bousquet}, 
B\'enyi et al. \cite{Benyi} showed that the  weak ascent sequences
can uniquely encode each of the following objects:
permutations avoiding a certain bivincular pattern of length $4$, upper triangular binary matrices satisfying a column restriction, and factorial posets that are ``special" $(3+1)$-free.
	
Very recently, 	Dukes and Sagan \cite{Dukes2023} introduced and studied  more general sequences which are called $d$-ascent sequences.
Given a sequence $x = x_1x_2\cdots x_n$ and an integer $d\geq 0$, an index $i$ ($1\leq i \leq n-1$) is said 
to be a {\em d-ascent} of $x$ if $x_{i+1} > x_i -d$. 
Let $\mathrm{dAsc}(x)$ denote the set of d-ascents of $x$ and let $\mathrm{dasc}(x)$ denote the number of d-ascents of $x$.
The sequence $x$ is called a {\em d-ascent sequence} if  $x_1=0$ and 
$0 \leq x_i \leq \mathrm{dasc}(x_1x_2\cdots x_{i-1})+1$ for all $2 \leq i \leq n$. 
	It is easily seen that  ascent sequences and weak ascent sequences correspond to the $d$-ascent sequences when $d = 0$ and $d = 1$, respectively.
	Let $\mathcal{A}^d_n$ denote the set of $d$-ascent sequences of length $n$.
	For example, we have  $x = 002143 \in \mathcal{A}^2_6$.
	It should be mentioned that $d$-ascent sequences are different from the 
	$p$-ascent sequences introduced by Kitaev and Remmel \cite{Kitaev2017}.

    Dukes and Sagan \cite{Dukes2023} showed that some maps from the weak case
    in \cite{Benyi} can be extended to bijections for general $d$ while the extensions of others continue to be injective but not surjective.
    To be specific, they constructed a bijection between $d$-ascent sequences and 
    upper triangular  matrices satisfying a column restriction and 
    a bijection between $d$-ascent sequences and matchings with restricted nestings.
    They also constructed an injection from $d$-ascent sequences to permutations
    avoiding a bivincular pattern of length $d+3$ and an injection from $d$-ascent sequences to factorial posets  avoiding a specially labeled poset with $d+3$ elements.
    
    The purpose of the present work is to complete the results of Dukes and Sagan 
    \cite{Dukes2023} by constructing a bijection between 
     a class of permutations
      which we call 
      difference $d$ permutations and $d$-ascent sequences (in Section \ref{sec:permutation}) and  a bijection between 
      a class of posets which we call  difference $d$ posets and $d$-ascent sequences
    (in Section \ref{sec:poset}).
    Our results are  extensions 
    of certain results of 
    Dukes and Sagan \cite{Dukes2023}.
	We also give an answer to a problem posed by Dukes and Sagan in the same paper by
	giving a direct bijection between two classes of matrices (in Section \ref{sec:matrix}).

\section{Permutations}\label{sec:permutation}

In this section, we will introduce a class of permutations that we call  
difference $d$ permutations and show that there is a bijection between difference $d$ permutations and $d$-ascent sequences.

Recall that $d$-ascent sequences are closely related to permutations avoiding a bivincular pattern.
The notion of pattern avoiding permutations was introduced by Knuth \cite{Knuth}
in 1970 to study the stack permutations.
Bousquet-M\'elou et al. \cite{Bousquet} initiated the study of bivincular patterns
and showed that ascent sequences are in bijection with permutations avoiding 
a particular bivincular pattern of length $3$.

For nonnegative integers $m,n$, we let $[m,n]=\{m,m+1,\ldots,n\}$, and when 
$m=1$ we abbreviate this to $[n]$.
Let $\mathcal{S}_n$ denote the set of permutations of 
$[n]$.
Given a permutation $\pi \in \mathcal{S}_n$ and a permutation $\sigma \in \mathcal{S}_k$,
an {\em occurrence} of $\sigma$ in $\pi$ is a subsequence $\pi_{i_1}\pi_{i_2}\cdots \pi_{i_k}$
of $\pi$ that is order isomorphic to $\sigma$.
We say $\pi$   {\em contains}  the (classical) pattern $\sigma$ if $\pi$ contains an occurrence of $\sigma$.
Otherwise, we say $\pi$ {\em avoids} the pattern $\sigma$  or
$\pi$ is {\em $\sigma$-avoiding}.
To contain a bivincular pattern $\sigma$, certain pairs of 
elements of the occurrence must be adjacent in $\pi$ and others must be 
adjacent as integers.
In the first case, we put a vertical bar between the elements of $\sigma$, 
and in the second case, we put a bar over the smaller of the two integers.
To illustrate, if $cdab$ is an occurrence  of the bivincular pattern $3| 41\bar{2}$ in $\pi$, then we have $a<b<c<d$ with $c$ and $d$ adjacent in $\pi$ and $c = b+1$.
For any pattern $\sigma$ (classical or bivincular),
let $\mathcal{S}_n(\sigma)$ denote the set of $\sigma$-avoiding permutations
of length $n$.
Define \[\tau_d = (d-1)|d12\cdots \overline{(d-2)}.\]

\begin{theorem}(\cite{Benyi,Bousquet}) \label{thm:pattern3-4}
	For $n \geq 1$, there is a bijection between $\mathcal{A}_n$ and  $\mathcal{S}_n(\tau_3)$, and  a bijection between $\mathcal{A}^1_n$ and  $\mathcal{S}_n(\tau_4)$.
\end{theorem}

\begin{theorem}(\cite{Dukes2023}, Theorem 4.5) \label{thm:pattern-general}
	For  $d\geq0$ and $n\geq 1$, there is an injection $\mathrm{pe}$ from $\mathcal{A}^d_n$ to
	$\mathcal{S}_n(\tau_{d+3})$.
\end{theorem}

The map pe induces a bijection between $\mathcal{A}_n$ and  $\mathcal{S}_n(\tau_3)$ when $d = 0$ and a bijection 
between $\mathcal{A}^1_n$ and  $\mathcal{S}_n(\tau_4)$
when $d = 1$.
Hence Theorem \ref{thm:pattern-general} is a generalization
of Theorem \ref{thm:pattern3-4}.
Dukes and Sagan \cite{Dukes2023} posed the following question.

\begin{question}(\cite{Dukes2023}, Question 8.4)\label{ques:1}
	Fix $d\geq 2$.
	Is there a set $\Sigma_d$ of bivincular patterns containing $\tau_{d}$ such that 
	$|\mathcal{A}^d_n| = |\mathcal{S}_n(\Sigma_{d+3})|$ for all $n\geq 1$?
\end{question}

This is actually the motivation and original intention behind our writing of this section.
However, instead of giving a direct answer to Question \ref{ques:1}, we 
introduce a class of permutations which we call difference $d$ permutations.
We show that difference $d$ permutations are $\tau_{d+3}$-avoiding
(in Theorem \ref{thm: d-permutation}) and in bijection with 
$d$-ascent sequences (in Theorem \ref{thm:phi}).

To introduce difference $d$ permutations, we will need the notion
of $d$-active elements in a permutation.
Let $d\geq 0, n\geq 1$, and let  $\pi= \pi_1\pi_2\cdots \pi_n$ be a permutation in $\mathcal{S}_n$.
We define the {\em $d$-active elements}
 of $\pi$ by the following procedure:
\begin{itemize}
	\item Set $1$ to be a $d$-active element of $\pi$.
	\item For $k = 2,3,\ldots,n$,
	   if $k$ is to the left of $k-1$ in $\pi$ and there exist at least $d$ $d$-active elements which are smaller than $k-1$ between $k$ and $k-1$ in $\pi$,   we say $k$ is a \emph{$d$-inactive element}  of $\pi$, otherwise, we say $k$ is  a $d$-active element of $\pi$.
\end{itemize}
In what follows, we abbreviate $d$-active (resp. $d$-inactive) to active (resp. inactive) if 
the value of $d$ is clear from the context.
Let $\mathrm{Act}(\pi)$ be the set of active elements of $\pi$ and 
let $\mathrm{act}(\pi)$ be the number of active elements of $\pi$.
For example, consider  $\pi = 42617385 \in \mathcal{S}_8$.
If $d = 0$, we have $\mathrm{Act}(\pi) = \{1,3,5,7,8\}$  and hence $\mathrm{act}(\pi) = 5$.
If $d = 2$, we have $\mathrm{Act}(\pi) = \{1,2,3,5,7,8\}$  and hence $\mathrm{act}(\pi) = 6$.

Given a permutation $\pi= \pi_1\pi_2\cdots \pi_n$, an index $i$ ($1 \leq i \leq n-1$) 
is said to be an {\em ascent} of $\pi$  if $\pi_{i+1} > \pi_i$, and we call $\pi_i$ an {\em ascent bottom} of $\pi$.
Let $\mathrm{Ascbot}(\pi)$ be the set of ascent bottoms of  $\pi$.
For example, let $\pi =  42617385$, we have $\mathrm{Ascbot}(\pi) = \{1,2,3\}$.

For $d\geq 0$, we call a permutation $\pi$ a {\em difference $d$ permutation}
if $\mathrm{Ascbot}(\pi) \subseteq \mathrm{Act}(\pi)$.
In what follows, we abbreviate difference $d$ permutations to difference permutations if 
the value of $d$ is clear from the context.
Let $\mathcal{S}^d_n$ denote the set of difference $d$ permutations in $\mathcal{S}_n$.
For example, let $\pi = 42617385$ be  the permutation  as given above.
If $d = 0$,  we have $\mathrm{Ascbot}(\pi) = \{1,2,3\} \nsubseteq  \{1,3,5,7,8\} = \mathrm{Act}(\pi)$.
Hence $\pi$ is not a  difference permutation in $\mathcal{S}^0_8$.
If $d = 2$,  we have $\mathrm{Ascbot}(\pi) = \{1,2,3\} \subseteq  \{1,2,3,5,7,8\} = \mathrm{Act}(\pi)$.
Hence $\pi$ is a difference permutation in $\mathcal{S}^2_8$.
It turns out that difference permutations are closely related to 
permutations avoiding bivincular patterns.

\begin{theorem}\label{thm: d-permutation}
	For $d\geq 0$ and $n \geq 1$, we have $\mathcal{S}^d_n \subseteq \mathcal{S}_n(\tau_{d+3})$
	and  equality holds when $d = 0$ or $d= 1$.
\end{theorem}

\pf
We first prove that $\mathcal{S}^d_n \subseteq \mathcal{S}_n(\tau_{d+3})$.
Let $\pi\in \mathcal{S}^d_n$.
We proceed to prove that $\pi  \in \mathcal{S}_n(\tau_{d+3})$.
If not, then there exists
an occurrence $ijk_1k_2\cdots k_d k_{d+1}$ of $\tau_{d+3}$ in $\pi$.
This means that $k_1<k_2<\cdots < k_d<k_{d+1}= i-1<i <j$ with $i,j$ adjacent in 
$\pi$.
Then we have $i \in \mathrm{Ascbot}(\pi)$.
For each $1\leq r \leq d$, 
since $k_r < k_{r+1} \leq i-1$,  there must be two adjacent elements $\pi _t<\pi_{t+1}$ with $\pi_t < i-1$ among the elements in the factor of $\pi$ from $k_r$ to $k_{r+1}$.
This implies that $\pi_t \in \mathrm{Ascbot}(\pi)$.
Since $\pi \in \mathcal{S}^d_n$, we have $\pi_t \in\mathrm{Ascbot}(\pi)\subseteq \mathrm{Act}(\pi)$.
It follows that $i$ is to the left of $i-1$ in $\pi$ and 
there are at least $d$ active elements between $i$ and $i-1$ in $\pi$ which are smaller than $i-1$. 
From the definition of inactive elements, $i$ is inactive.
Since $i \in \mathrm{Ascbot}(\pi)$,  we have $\mathrm{Ascbot}(\pi) \nsubseteq
\mathrm{Act}(\pi)$, a contradiction.
Thus $\mathcal{S}^d_n \subseteq \mathcal{S}_n(\tau_{d+3})$.

It remains to prove that $\mathcal{S}_n(\tau_{d+3}) \subseteq \mathcal{S}^d_n$ for 
$d = 0,1$.
We only consider the case $d= 1$ as the other case $d = 0$ can be proven similarly.
Let $\pi\in \mathcal{S}_n(\tau_{4})$.
If $\pi \nin \mathcal{S}^1_n$, then there is some $\pi_k \in \mathrm{Ascbot}(\pi)$ such that $\pi_k \nin \mathrm{Act}(\pi)$.
This implies that $\pi_k < \pi_{k+1}$ and $\pi_k$ is to the left of $\pi_k-1$ with at least one active element $\pi_{\ell} < {\pi_k-1}$  between them.
Then $\pi_k\pi_{k+1} \pi_{\ell} (\pi_k-1)$ forms an occurrence of $\tau_4$ in $\pi$,
a contradiction to the fact that $\pi\in \mathcal{S}_n(\tau_{4})$.
This completes the proof.
\qed

For $d\geq 2$ and $n \geq d+3$, 
we remark that  $\mathcal{S}^d_n \subsetneq \mathcal{S}_n(\tau_{d+3})$  since one can check that 
the permutation  $\pi =(d+2)(d+3)\cdots n d\cdots 21(d+1)$
is a permutation in  $\mathcal{S}_n(\tau_{d+3})$ but not in 
$\mathcal{S}^d_n$.

\begin{lemma}\label{lem:perm_gene}
	Given $d\geq 0$ and $n\geq 2$, let $\sigma$ be a permutation in $\mathcal{S}_{n-1}$ and let $\pi$ be a permutation obtained from $\sigma$ by inserting the element $n$ into $\sigma$.
	Then $\pi \in \mathcal{S}^d_n$ if and only if $\sigma\in \mathcal{S}^{d}_{n-1}$
	and $n$ is inserted  before $\sigma$ or after some active element of $\sigma$.
\end{lemma}

\pf
Suppose  that $\pi \in \mathcal{S}^d_n$. 
By the definition of difference $d$ permutations, we have $\mathrm{Ascbot}(\pi) \subseteq \mathrm{Act}(\pi)$.
Notice that the elements of $\sigma$ remain active or 
inactive after
the insertion of $n$ into $\sigma$.
If the element $n$ is inserted after some inactive  element
$j$ of $\sigma$, then $j\in \mathrm{Ascbot}(\pi)$ but 
$j \nin \mathrm{Act}(\pi)$, a contradiction to the fact that 
$\mathrm{Ascbot}(\pi) \subseteq \mathrm{Act}(\pi)$.
Hence $n$ is inserted  before $\sigma$ or after some active element of $\sigma$.
We now show  that $\sigma \in \mathcal{S}^{d}_{n-1}$.
If not, there is some $k \in \mathrm{Ascbot}(\sigma)$ but $k \nin\mathrm{Act}(\sigma)$.
Hence $k\nin \mathrm{Act}(\pi)$.
It is easily seen that $\mathrm{Ascbot}(\sigma) \subseteq \mathrm{Ascbot}(\pi)$,
thereby $k \in \mathrm{Ascbot}(\pi)$.
Then we have $\mathrm{Ascbot}(\pi) \nsubseteq \mathrm{Act}(\pi)$,
a contradiction.
Hence $\sigma \in \mathcal{S}^{d}_{n-1}$.

For the converse,  suppose that $\sigma\in \mathcal{S}^{d}_{n-1}$.
We have two cases.
If $n$ is inserted  before $\sigma$,
then we have $\mathrm{Ascbot}(\pi) = \mathrm{Ascbot}(\sigma) \subseteq \mathrm{Act}(\sigma) \subseteq  \mathrm{Act}(\pi)$.
Thus $\pi$ is a difference $d$ permutation.
If  $n$ is inserted  after some active element $i$ of $\sigma$,
then the newly (possibly) added  ascent bottom $i$  is an active element of $\pi$.
From the fact that $\mathrm{Ascbot}(\sigma) \subseteq \mathrm{Act}(\sigma)$ we  also have  $\mathrm{Ascbot}(\pi) \subseteq \mathrm{Act}(\pi)$, 
namely, $\pi$ is a difference $d$ permutation.
This completes the proof.
\qed

Based on Lemma \ref{lem:perm_gene}, 
we now define a map $\phi$ from  difference $d$ permutations $\mathcal{S}^d_n$ to $d$-ascent sequences $\mathcal{A}^d_n$.
Our map $\phi$ is defined recursively.
For $n = 1$, we define $\phi(1) = 0$. 
Next let $n \geq 2$ and suppose that $\pi$ is obtained from $\sigma$ by inserting the element $n$ after the $x_n$-th active element of $\sigma$ (reading from left to right).
We set $x_n = 0$ if $n$ is inserted before $\sigma$.
Then the sequence associated to $\pi$ is $\phi(\pi) = x_1x_2\cdots x_n$,
where $x_1x_2\cdots x_{n-1} = \phi(\sigma)$.

\begin{example}
	Let $d=2$ and let $\pi = 42617385$ be a difference  permutation in $\mathcal{S}^2_8$.
	Then we have $\phi(\pi) = 00203124$ with the following recursive insertion 
	of new maximal values.
	The elements colored by red indicate the active elements.
	 \begin{align*}
		\red{1}
		&\,\xrightarrow{x_2=0}\,  \red{2} \red{1} \\
		&\,\xrightarrow{x_3=2}\,  \red{2} \red{1} \red{3} \\
		&\,\xrightarrow{x_4=0}\,  4\red{2} \red{1} \red{3}  \\
		&\,\xrightarrow{x_5=3}\,  4\red{2} \red{1} \red{3} \red{5} \\
		&\,\xrightarrow{x_6=1}\,  4\red{2} 6\red{1} \red{3} \red{5}  \\
		&\,\xrightarrow{x_7=2}\,  4\red{2} 6\red{1} \red{7} \red{3} \red{5} \\
		&\,\xrightarrow{x_8=4}\,  4\red{2} 6\red{1} \red{7} \red{3} \red{8}\red{5}  .
	\end{align*}
\end{example}

From the construction of the map $\phi$, it is easily seen that 
$x_i$ $(1\leq i \leq n)$ is the number of active elements to the left of $i$ in $\pi$ which are smaller than $i$.

\begin{theorem}\label{thm:phi}
	For $d\geq 0$ and $n\geq 1$, 
	the map $\phi$ is a bijection between $\mathcal{S}^d_n$ and $\mathcal{A}^d_n$.
	Furthermore, we have $\mathrm{act}(\pi) = \mathrm{dasc}(\phi(\pi)) + 1$ for all $\pi \in \mathcal{S}^d_n$.
\end{theorem}

\pf
Since the sequence $\phi(\pi)$ encodes the construction of $\pi$, 
the map $\phi$ is injective.
To prove $\phi$ is a bijection,  we need to show that
the image  $\phi(\mathcal{S}^d_n)$ is the set $\mathcal{A}^d_n$.
The recursive construction of the map $\phi$ tells us that 
$x = x_1x_2\cdots x_n\in \phi(\mathcal{S}^d_n)$  if and only if 
$x' = x_1x_2\cdots x_{n-1} \in \phi(\mathcal{S}^d_{n-1})$ and 
$0\leq x_n \leq \mathrm{act}(\phi^{-1}(x'))$.
By induction on $n$ and the definition of $d$-ascent sequences, 
to  prove $\phi(\mathcal{S}^d_n) = \mathcal{A}^d_n$, 
it is sufficient to show that $\mathrm{act}(\pi) = \mathrm{dasc}(\phi(\pi)) + 1$.

Let us focus on the property $\mathrm{act}(\pi) = \mathrm{dasc}(\phi(\pi)) + 1$.
We will prove the result by induction on $n$ where $n = 1$ is trivial.
Assume the result for $n-1$.
We need to prove the result for $n$.
Let $\pi = \pi_1\pi_2\cdots \pi_n$ be a permutation in $\mathcal{S}^d_n$ and 
$x = x_1x_2\cdots x_n = \phi(\pi)$.
Then $x' = x_1x_2\cdots x_{n-1} = \phi(\sigma)$, where 
$\sigma$ is the permutation obtained from $\pi$ by deleting the element 
$n$ from $\pi$. 
By the induction hypothesis,  we have $\mathrm{act}(\sigma) = \mathrm{dasc}(x') + 1$.
To prove $\mathrm{act}(\pi) = \mathrm{dasc}(\phi(\pi)) + 1$,
it suffices to show that $n\nin \mathrm{Act}(\pi)$ if and only if
$n-1 \nin \mathrm{dAsc}(x)$.
Recall that $x_i$ $(1\leq i \leq n)$ is the number of active elements to the left of $i$ in $\pi$ which are smaller than $i$.
If $n$ is to the left of $n-1$ in $\pi$, we have 
 that $x_{n-1}$ is the sum of $x_n$ and the number of active elements between $n$ and $n-1$ in $\pi$ which are smaller than $n-1$.
 Then by the definition of inactive elements, $n \nin \mathrm{Act}(\pi)$ 
if and only if  $x_n \leq x_{n-1} -d$, namely $n-1 \nin \mathrm{dAsc}(x)$.
This completes the proof.
\qed

Combining Theorems \ref{thm: d-permutation} and \ref{thm:phi} 
gives  new proofs of Theorems \ref{thm:pattern3-4} and \ref{thm:pattern-general}.
Let $\pi$ be a permutation in $\mathcal{S}_n$ and let  $\sigma$ be 
a bivincular pattern.
Fix $d \geq 0$.
We generalize the concept of bivincular patterns  by further specifying that certain elements of an occurrence of the pattern $\sigma$ must be $d$-active in $\pi$.
We indicate this by underlining the corresponding elements in $\sigma$.
For instance, if $\pi_{i_1}\pi_{i_2}\pi_{i_3}\pi_{i_4}\pi_{i_5}$ is an occurrence of the special bivincular pattern 
$4| 5\underline{1}\underline{2}\bar{3}$ in $\pi$, then we have $\pi_{i_3}<\pi_{i_4}<\pi_{i_5}<\pi_{i_1}<\pi_{i_2}$ with $\pi_{i_1}$ and $\pi_{i_2}$ adjacent in $\pi$,  $\pi_{i_1} = \pi_{i_5}+1$, and $\pi_{i_3},\pi_{i_4}$  are $d$-active elements in $\pi$.
Let 
\[
\Sigma_{d+3} = \left\{ \tau'_{\mu} \mid \mu \in \mathcal{S}_d\right\},
\]
where $$\tau'_{\mu} = (d+2)|(d+3)\underline{\mu} \overline{d+1}.$$
Let $\mathcal{S}_n(\Sigma_{d+3})$ be the set of permutations that avoid all the special bivincular patterns in $\Sigma_{d+3}$. 
By the definition of difference $d$ permutations,
it is easily seen that $\mathcal{S}^d_n = \mathcal{S}_n(\Sigma_{d+3})$.
Hence our results can be regard as an answer to Question \ref{ques:1}.
It should be mentioned that the map $\phi^{-1}$ is different from 
the map pe in Theorem \ref{thm:pattern-general}.


\section{Posets}\label{sec:poset}

In this section, 
we will introduce a class of posets which we call difference $d$ posets and 
show that there is a bijection between difference $d$ posets and $d$-ascent sequences.

Let $P$ be a poset (partial ordered set).
We say $P$ is {\em $(a+b)$-free} if it does not contain an 
(induced) subposet which is isomorphic to the disjoint union of an $a$-element 
chain and a $b$-element chain.
For example, the poset $P$ whose Hasse diagram is on the left of Figure \ref{fig:posets} is not $(3+1)$-free because the subposet of $P$ consisting of the elements
$\{3,5,7,8\}$ forms an occurrence of $(3+1)$.
But the poset  $P$  is $(2+2)$-free.

\begin{figure}[H]
	\begin{center}
		\begin{tikzpicture}[font =\small , scale = 0.45, line width = 0.7pt]
			\foreach \i / \j \k in {1/1/1,3/1/2,5/1/4,1/3/3,3/3/6,5/3/7,3/5/5,5/7/8} {\draw[black,fill =pink](\i,\j)circle(10pt);
				\node at(\i,\j){\k};
			}
			\tikzmath{\c = 0.38;};
			\draw (1,1+\c)--(1,3-\c);
			\draw (1,1+\c)--(3,3-\c);
			\draw (1,1+\c)--(5,3-\c);
			\draw (3,1+\c)--(1,3-\c);
			\draw (3,1+\c)--(5,3-\c);
			\draw (5,1+\c)--(3,5-\c);
			\draw (1,3+\c)--(3-0.2,5-\c+0.1);
			\draw (3,3+\c)--(5,7-\c);
			\draw (3+0.2,5+\c-0.1)--(5,7-\c);
			\node at (3,-1) {$P$};
			
			\tikzmath{\x = 10;};
			\foreach \i \j \k in 
			{1/1/1,3/1/2,5/1/5,2/3/3,4/3/4,3/5/6}
			{\draw[black,fill =pink](\i+\x,\j)circle(10pt);
				\node at(\i+\x,\j){\k};
			}
			\draw (1+\x,1+\c)--(2+\x-0.2,3-\c+0.1);
			\draw (1+\x,1+\c)--(4+\x,3-\c);
			\draw (3+\x,1+\c)--(4+\x,3-\c);
			\draw (2+\x+0.1,3+\c)--(3+\x,5-\c);
			\draw (4+\x,3+\c)--(3+\x,5-\c);
			\node at (3+\x,-1) {$Q$};
		\end{tikzpicture}
	\end{center}
	\caption{Two factorial posets.} \label{fig:posets}
\end{figure}

Let $P$ be a poset on integers. 
We will use $<_P$ to denote the partial order on $P$ and $<$ for the total order on the integers.
We call $P$ {\em compatible} if it satisfies the following rule:
\[i<_Pj \Longrightarrow i<j\]
for all  $i,j \in P$.
We call a  poset $P$ on $[n]$ a {\em factorial poset} if it satisfies the following rule:
\begin{align}\label{factorial-rule}
	i<j \text{ and } j<_Pk \; \Longrightarrow \; i<_P k
\end{align}
for all $i,j,k\in [n]$.
Factorial posets were first introduced by Claesson and Linusson \cite{Claesson2} and are easily seen to be compatible and $(2+2)$-free.
The reason they are called factorial posets is that there is a natural bijection $\omega$ between factorial posets and inversion sequences.
Given a factorial poset $P$ on $[n]$, define $\omega(P) = a_1a_2\cdots a_n$, where
\begin{align*}
a_i = 
\begin{cases}
 		 	0,  &\text{if $i$ is a minimal element of } P, \\
 		\text{max}\{j \mid j <_P i\},  &\text{otherwise.}
\end{cases}
\end{align*}
Let $\mathcal{A}(P)$ be the set of
 nonzero elements of $\omega(P)$.
That is, $\mathcal{A}(P) = \{a_i \in \omega(P)\mid a_i > 0\}$.
For example, the two posets shown in Figure \ref{fig:posets} are both factorial 
posets.
Moreover, we have $\omega(P) = 00204126$ and $\omega(Q) = 001204$,
from which it follows that $\mathcal{A}(P) = \{1,2,4,6\}$ and $\mathcal{A}(Q) = \{1,2,4\}$.

B\'enyi et al. \cite{Benyi} built a bijection between  weak ascent sequences
and factorial posets which do not contain a specially labeled $(3+1)$ subposet.
In order to extend this result, 
Dukes and Sagan \cite{Dukes2023} introduced a  special
compatible poset $P_d$ which is the disjoint union of a $(d-1)$-element chain
 and an isolated element whose label is one more than the second largest element of the $(d-1)$-element chain. 
Let $P$ be a factorial poset. 
We say $P$ {\em contains} the special poset $P_d$ if there exists a 
subposet of $P$ which is the disjoint union of a $(d-1)$-element chain
$$i_1<_P i_2 <_P \cdots <_P i_{d-2} <_P i_{d-1}$$
and an isolated element $i_{d-2}+1$.
Otherwise, we say $P$ is {\em special $P_d$-free}.
Let $\mathcal{P}_n(P_d)$ denote the set of special $P_d$-free 
factorial posets on $[n]$.
For example, one can check (carefully) that the poset $P$ shown in Figure \ref{fig:posets} 
is special $P_4$-free.
And the poset $Q$ shown in Figure \ref{fig:posets}
contains the special poset $P_4$ because the subposet
of $Q$ consisting of  the elements $\{1,4,5,6\}$ forms  an occurrence of the special poset $P_4$.
The following theorem is an extension of the map of B\'enyi et al.
\cite{Benyi} from the weak case.

\begin{theorem}(\cite{Dukes2023}, Theorem 5.4)\label{thm:poset-injection}
     For $d \geq 1$ and $n \geq 1$, there is an injection $\mathrm{po}$ from 
     $\mathcal{A}^d_n$ to $\mathcal{P}_n(P_{d+3})$.
\end{theorem}

In analogy to permutations, 
Dukes and Sagan \cite{Dukes2023} posed the following question.


\begin{question}(\cite{Dukes2023}, Question 8.4)\label{ques:2}
	Fix $d\geq 2$.
	Is there a set $\Sigma'_d$ of special posets containing $P_{d}$ such that 
	$|\mathcal{A}^d_n| = |\mathcal{P}_n(\Sigma'_{d+3})|$ for all $n\geq 1$?
\end{question}

The purpose of this section is to give an answer to 
Question \ref{ques:2}.
To this end, we introduce a class of posets which 
we call  difference $d$ posets.
We show that difference $d$ posets are special $P_{d+3}$-free
(in Theorem \ref{thm:d-poset}) and in bijection with 
$d$-ascent sequences (in Theorem \ref{thm:psi}).

Analogously to permutations, we need to define
$d$-active (active)  elements on factorial posets.
Let $d\geq 0$ and let $P$ be a factorial poset on $[n]$ with $
$$\omega(P) = a_1a_2\cdots a_n$.
We define the {\em active elements}  
of $P$ in the following procedure:
For $k=1,2,\ldots n-1$,  if $a_{k+1}\leq a_k$ and there are at least $d$ active
elements in $[a_{k+1}+1, a_k]$, we say $k$ is an \emph{inactive element}  of $P$.
Equivalently, if the set $\{u \in P\mid u<_P k, u\nless_P k+1\}$ contains at least $d$ active elements, 
we say $k$ is inactive.
Otherwise, we say $k$ is an active element of $P$.
Set $n$ to be an inactive element of $P$.
Let $\mathrm{Act}(P)$ denote the set of active elements of $P$ and 
let $\mathrm{act}(\pi)$ denote the number of active elements of $P$.
For example,  let $P$ be the poset shown in Figure \ref{fig:posets} 
with $\omega(P) = 00204126$.
If $d = 0$, we have $\mathrm{Act}(P) = \{2,4,6,7\}$ and hence 
$\mathrm{act}(P) = 4$.
If $d = 2$, we have $\mathrm{Act}(P) = \{1,2,4,6,7\}$ and hence 
$\mathrm{act}(P) = 5$.

We are now in a position to define difference $d$ posets.
For  $d\geq 0$,  we call a factorial poset $P$ a
{\em difference $d$ poset} if $\mathcal{A}(P) \subseteq \mathrm{Act}(P)$.
In what follows, we abbreviate difference $d$ posets to difference posets if 
the value of $d$ is clear from the context.
Let $\mathcal{P}^d_n$ denote the set of difference $d$ posets on $[n]$.
For example, let $d=2$ and let $P,Q$ be the two posets shown in Figure \ref{fig:posets}.
For poset $P$, we have $\mathcal{A}(P) = \{1,2,4,6\} \subseteq  \{1,2,4,6,7\} = \mathrm{Act}(P)$.
Hence $P$ is  a  difference  poset in $\mathcal{P}^2_8$.
For poset $Q$, we have $\mathcal{A}(Q) = \{1,2,4\} \nsubseteq  \{1,2,3,5\} = \mathrm{Act}(Q)$.
Hence $Q$ is not a difference poset in $\mathcal{P}^2_6$.
It turns out that difference $d$ posets are closely related to special 
$P_d$-free posets.

\begin{theorem}\label{thm:d-poset}
	For all $d,n \geq 1$, we have $\mathcal{P}^d_n \subseteq \mathcal{P}_n(P_{d+3})$
	and equality holds when  $d = 1$.
\end{theorem}

\pf 
We first show that $\mathcal{P}^d_n \subseteq \mathcal{P}_n(P_{d+3})$.
Let $P$ be a factorial poset in $\mathcal{P}^d_n$ with $\omega(P) = a_1a_2\cdots a_n$.
By the definition of difference $d$ posets, we have $\mathcal{A}(P)\subseteq \mathrm{Act}(P)$.
We need to prove  $P \in \mathcal{P}_n(P_{d+3})$.
If not, then there is 
 a subposet  of $P$ which is the disjoint union of a $(d+2)$-element chain
$$i_1<_P i_2 <_P \cdots <_P i_{d+1} <_P i_{d+2}$$
and an isolated element $i_{d+1}+1$.
By the definition of $\omega(P)$ and  rule (\ref{factorial-rule}), 
it is easily seen that $$i_1\leq a_{i_2}<i_2\leq a_{i_3}<\cdots <i_d \leq a_{i_{d+1}}
<i_{d+1}.$$
Since $i_{d+1} <_P  i_{d+2}$ and $i_{d+1} +1  \nless _P  i_{d+2}$,
we have $a_{i_{d+2}} = i_{d+1}$.
Combining the fact that
$\mathcal{A}(P)\subseteq \mathrm{Act}(P)$, 
we obtain that $a_{i_k} \in \mathrm{Act}(P)$ for $2\leq k\leq d+2$.
In particular, we have $i_{d+1}=a_{i_{d+2}}\in \mathrm{Act}(P)$.
Note that $i_1 \nless_P i_{d+1}+1$. 
By rule (\ref{factorial-rule}),
we have $a_{i_k} <_P i_{d+1}$ and $a_{i_k}\nless_P i_{d+1}+1$ for $2\leq k \leq d+1$.
Then by the definition of inactive elements,  $i_{d+1}$ is  an inactive element of $P$, a contradiction.
Thus, $\mathcal{P}^d_n \subseteq \mathcal{P}_n(P_{d+3})$.

It remains to prove that $\mathcal{P}_n(P_{4}) \subseteq  \mathcal{P}^1_n$.
Let $P\in \mathcal{P}_n(P_{4})$. 
If $P \nin \mathcal{P}^1_n$, then there is some $i$ such that $i\in \mathcal{A}(P)$
but $i \nin \mathrm{Act}(P)$.
From the definition of $\mathcal{A}(P)$, there is some $k\in [n]$ such that $a_k = i$.
By the definition of $\omega(P)$,  we have 
$i <_P k$ and $i+1 \nless_P k$.
From the definition of inactive elements, 
there is some active element $j$ satisfying that $j<_P i$ but $j\nless_P i+1$.
It is easily seen that $i+1\neq k$.
Then the subposet of $P$ consisting of the elements $\{j,i,i+1,k\}$ forms an 
occurrence of the special poset $P_4$,  a contradiction.
This completes the proof.
\qed

Numerical evidence shows that $\mathcal{P}^0_n \nsubseteq \mathcal{P}_n(P_{3})$.
For example, let $P$ be the poset shown in Figure \ref{fig:couterexample}.
One can check that $P\in \mathcal{P}^0_5$ but $P \nin \mathcal{P}_5(P_{3})$.
However, we have the following theorem.

\begin{figure}[H]
	\begin{center}
		\begin{tikzpicture}[font =\small , scale = 1, line width = 0.7pt]
			\foreach \i / \j \k in {2/0/1,2/1/2,1/1/3,3/1/4,2/2/5} {\draw[black,fill =pink](\i,\j)circle(5pt);
				\node at(\i,\j){\k};
			}
			\tikzmath{\c = 0.84; \d = 0.87;};
			\draw (2,0+\c)--(2,1-\c);
			\draw (2,1+\c)--(2,2-\c);
			\draw (1+\d,1+\d)--(2-\d,2-\d);
			\draw (2+\d,\d)--(3-\d,1-\d);
		\end{tikzpicture}
	\end{center}
	\caption{A poset in $\mathcal{P}^0_5$ but not in $\mathcal{P}_5(P_{3})$.
} \label{fig:couterexample}
\end{figure}

\begin{theorem}
	For $n \geq 1$, we have $\mathcal{P}_n(P_{3}) \subseteq \mathcal{P}^0_n$.
\end{theorem}

\pf 
Let $P\in \mathcal{P}_n(P_{3})$ with $\omega(P) = a_1a_2\cdots a_n$. 
If $P \nin \mathcal{P}^0_n$, then there is some $i$ such that $i\in \mathcal{A}(P)$
but $i \nin \mathrm{Act}(P)$.
By the definition of $\mathcal{A}(P)$,  we have 
$i <_P k$ and $i+1 \nless_P k$ for some $k$.
From the definition of inactive elements, we have $a_{i+1} \leq a_i$.
It follows that $i \nless_P i+1$.
It is easily seen that $i+1\neq k$.
Then the subposet of $P$ consisting of the elements $\{i,i+1,k\}$ forms an 
occurrence of the special poset $P_3$,  a contradiction.
This completes the proof.
\qed

Given a factorial poset $P$ on $[n]$ and $1\leq i \leq n$, 
let $P[i]$ denote the subposet of $P$ consisting of all the elements in $[i]$.
It is easily seen that $P[i]$ is also a 
factorial poset for $1\leq i\leq n$.
We proceed now to construct a bijection between difference $d$ posets and 
$d$-ascent sequences. 
The following lemma is needed.

\begin{lemma}\label{lem:poset_gene}
	Given $d\geq 0$ and $n\geq 2$, let $P$ be a factorial poset on $[n]$ with 
	$\omega(P) = a_1a_2 \cdots a_n$.
	Then $P\in \mathcal{P}^d_n$ if and only if $P[n-1] \in \mathcal{P}^d_{n-1}$ 
	and $a_n \in \{0,n-1\} \cup \mathrm{Act}(P[n-1])$.
\end{lemma}

\pf
Suppose that $P\in \mathcal{P}^d_n$.
By the definition of difference $d$ posets, we have $\mathcal{A}(P)\subseteq \mathrm{Act}(P)$.
If $a_n \neq 0$ and  $a_n \neq n-1$, then $a_n \in \mathcal{A}(P) \subseteq \mathrm{Act}(P)$.
Notice that  either $\mathrm{Act}(P) = \mathrm{Act}(P[n-1])$ or
$\mathrm{Act}(P) = \mathrm{Act}(P[n-1]) \cup \{n-1\}$.
Hence we have $a_n \in \mathrm{Act}(P[n-1])$.
We now show that $P[n-1] \in \mathcal{P}^d_{n-1}$.
By the definition of $\mathcal{A}(P)$, we have 
 $\mathcal{A}(P[n-1]) \subseteq \mathcal{A}(P) \subseteq \mathrm{Act}(P)$.
Since all the elements in $\mathcal{A}(P[n-1])$ are smaller than $n-1$, 
we derive that  $\mathcal{A}(P[n-1])  \subseteq \mathrm{Act}(P[n-1])$.
Thus $P[n-1] \in \mathcal{P}^d_{n-1}$.

For the converse, suppose that $P[n-1] \in \mathcal{P}^d_{n-1}$ and $a_n \in \{0,n-1\} \cup \mathrm{Act}(P[n-1])$.
Then we have $\mathcal{A}(P[n-1])\subseteq \mathrm{Act}(P[n-1])$.
We need to prove $P\in \mathcal{P}^d_n$, namely, 
$\mathcal{A}(P) \subseteq \mathrm{Act}(P)$.
There are three cases.
If $a_n= 0$,  then we have $\mathcal{A}(P) = \mathcal{A}(P[n-1]) \subseteq \mathrm{Act}(P[n-1]) \subseteq \mathrm{Act}(P)$.
If $a_n = n-1$, then we have $\mathcal{A}(P) = \mathcal{A}(P[n-1]) \cup \{n-1\} \subseteq \mathrm{Act}(P[n-1]) \cup \{n-1\} \subseteq \mathrm{Act}(P) \cup \{n-1\}$.
Since $a_n=n-1 > a_{n-1}$, we deduce that $n-1$ is an active element of $P$.
Thus, we have $\mathcal{A}(P)\subseteq \mathrm{Act}(P)$.
If $a_n \in \mathrm{Act}(P[n-1])$,  it is easily seen that 
$\mathcal{A}(P) = \mathcal{A}(P[n-1]) \cup \{a_n\} \subseteq \mathrm{Act}(P[n-1]) \subseteq \mathrm{Act}(P)$.
This completes the proof. 
\qed

Lemma \ref{lem:poset_gene} enables us to construct a map $\psi$ from difference $d$ posets to $d$-ascent sequences recursively.
Let $d\geq 0$ and let $P$ be a difference $d$ poset of $\mathcal{P}^d_n$ with  
$\omega(P) = a_1a_2\cdots a_n$.
For $n=1$, set $\psi(P) = 0$.
For $n \geq 2$, 
suppose that we have obtained $\psi(P[n-1]) = x_1x_2\cdots x_{n-1}$.
Then define $\psi(P) = x_1\cdots x_{n-1}x_n$,
where
\begin{align}\label{equ:psi}
x_n =
\begin{cases}
	0,  & \text{if } a_n = 0,\\
	\mathrm{act}(P[n-1])+1, &\text{if }a_n = n-1,\\
	i,   &\text{if }a_n \text{ is the $i$-th smallest element in } \mathrm{Act}(P[n-1]).
\end{cases}
\end{align}

\begin{example}\label{exam:poset}
	Let $d = 2$ and let $P$ be the difference $d$ poset shown in Figure \ref{fig:posets}.
	Then we have $\psi(P) = 00203124$ with the recursive construction in Figure \ref{fig:psi} in which the active elements are colored by green in each poset.
\end{example}

\begin{figure}[H]
	\begin{center}
		\begin{tikzpicture}[font =\small , scale = 0.47, line width = 0.7pt]
			\tikzmath{\x = -3; \y = -8;\c = 0.38;};
			\draw[black,fill = pink](1+\x,1)circle(10pt);
			\node at(1+\x,1){1};
			\tikzmath{\x = 3;};
			\foreach \i / \j \k in {1/1/1,3/1/2}
			 {	\ifthenelse{\k = 1} 
			 	{\draw[black,fill = green](\i+\x,\j)circle(10pt);}
			 	{\draw[black,fill = pink](\i+\x,\j)circle(10pt);}
					\node at(\i+\x,\j){\k};
			}
		    \tikzmath{\x = 11;};
   			\foreach \i / \j \k in {1/1/1,3/1/2}
		    {	
		    	\draw[black,fill = green](\i+\x,\j)circle(10pt);
		    	\node at(\i+\x,\j){\k};
		    }
    		\foreach \i / \j \k in {1/3/3}
		    {	
		    	\draw[black,fill = pink](\i+\x,\j)circle(10pt);
		    	\node at(\i+\x,\j){\k};
		    }
   			\draw(1+\x,1+\c)--(1+\x,3-\c);
	        \draw(3+\x,1+\c)--(1+\x,3-\c);
			\tikzmath{\x = 18;};
			\foreach \i / \j \k in {1/1/1,3/1/2}
			{	
				\draw[black,fill = green](\i+\x,\j)circle(10pt);
				\node at(\i+\x,\j){\k};
			}
			\foreach \i / \j \k in {1/3/3,5/1/4}
			{	
				\draw[black,fill = pink](\i+\x,\j)circle(10pt);
				\node at(\i+\x,\j){\k};
			}
   			\draw(1+\x,1+\c)--(1+\x,3-\c);
			\draw(3+\x,1+\c)--(1+\x,3-\c);
			
			\tikzmath{\x = 18;};
			\foreach \i / \j \k in {1/1/1,3/1/2,5/1/4}
			{	
				\draw[black,fill = green](\i+\x,\j+\y)circle(10pt);
				\node at(\i+\x,\j+\y){\k};
			}
			\foreach \i / \j \k in {1/3/3,3/5/5}
			{	
				\draw[black,fill = pink](\i+\x,\j+\y)circle(10pt);
				\node at(\i+\x,\j+\y){\k};
			}
			\draw(1+\x,1+\c+\y)--(1+\x,3-\c+\y);
			\draw(3+\x,1+\c+\y)--(1+\x,3-\c+\y);
			\draw (1+\x,3+\c+\y)--(3-0.2+\x,5-\c+0.1+\y);
			\draw (5+\x,1+\c+\y)--(3+\x,5-\c+\y);
			
			\tikzmath{\x = 11;};
			\foreach \i / \j \k in {1/1/1,3/1/2,5/1/4}
			{	
				\draw[black,fill = green](\i+\x,\j+\y)circle(10pt);
				\node at(\i+\x,\j+\y){\k};
			}
			\foreach \i / \j \k in {1/3/3,3/5/5,3/3/6}
			{	
				\draw[black,fill = pink](\i+\x,\j+\y)circle(10pt);
				\node at(\i+\x,\j+\y){\k};
			}
			\draw(1+\x,1+\c+\y)--(1+\x,3-\c+\y);
			\draw(3+\x,1+\c+\y)--(1+\x,3-\c+\y);
			\draw (1+\x,3+\c+\y)--(3-0.2+\x,5-\c+0.1+\y);
			\draw (5+\x,1+\c+\y)--(3+\x,5-\c+\y);
			\draw (1+\x,1+\c+\y)--(3+\x,3-\c+\y);
			
			\tikzmath{\x = 3;};
			\foreach \i / \j \k in {1/1/1,3/1/2,5/1/4,3/3/6}
			{	
				\draw[black,fill = green](\i+\x,\j+\y)circle(10pt);
				\node at(\i+\x,\j+\y){\k};
			}
			\foreach \i / \j \k in {1/3/3,3/5/5,5/3/7}
			{	
				\draw[black,fill = pink](\i+\x,\j+\y)circle(10pt);
				\node at(\i+\x,\j+\y){\k};
			}
			\draw(1+\x,1+\c+\y)--(1+\x,3-\c+\y);
			\draw(3+\x,1+\c+\y)--(1+\x,3-\c+\y);
			\draw (1+\x,3+\c+\y)--(3-0.2+\x,5-\c+0.1+\y);
			\draw (5+\x,1+\c+\y)--(3+\x,5-\c+\y);
			\draw (1+\x,1+\c+\y)--(3+\x,3-\c+\y);
			\draw (1+\x,1+\c+\y)--(5+\x,3-\c+\y);
			\draw (3+\x,1+\c+\y)--(5+\x,3-\c+\y);
			
			\tikzmath{\x = -5;};
			\foreach \i / \j \k in {1/1/1,3/1/2,5/1/4,3/3/6,5/3/7}
			{	
				\draw[black,fill = green](\i+\x,\j+\y)circle(10pt);
				\node at(\i+\x,\j+\y){\k};
			}
			\foreach \i / \j \k in {1/3/3,3/5/5,5/7/8}
			{	
				\draw[black,fill = pink](\i+\x,\j+\y)circle(10pt);
				\node at(\i+\x,\j+\y){\k};
			}
			\draw(1+\x,1+\c+\y)--(1+\x,3-\c+\y);
			\draw(3+\x,1+\c+\y)--(1+\x,3-\c+\y);
			\draw (1+\x,3+\c+\y)--(3-0.2+\x,5-\c+0.1+\y);
			\draw (5+\x,1+\c+\y)--(3+\x,5-\c+\y);
			\draw (1+\x,1+\c+\y)--(3+\x,3-\c+\y);
			\draw (1+\x,1+\c+\y)--(5+\x,3-\c+\y);
			\draw (3+\x,1+\c+\y)--(5+\x,3-\c+\y);
			\draw (3+\x,3+\c+\y)--(5+\x,7-\c+\y);
			\draw (3+0.2+\x,5+\c-0.1+\y)--(5+\x,7-\c+\y);
			
			\draw[->] (-0.5,1)--(2,1) node[font = \small,left = 6mm,above]{$x_2 = 0$};
			\tikzmath{\x =8;};
			\draw[->] (-0.5+\x,1)--(2+\x,1) node[font = \small,left = 6mm,above]{$x_3 = 2$};
			\tikzmath{\x =16;};
			\draw[->] (-0.5+\x,1)--(2+\x,1) node[font = \small,left = 6mm,above]{$x_4 = 0$};
			\tikzmath{\x =19;};
			\draw[->] (2+\x,0)--(2+\x,-2) node[font = \small,above = 5mm,right]{$x_5 = 3$};
			\tikzmath{\x =16.2;\y = -6;};
			\draw[->] (2+\x,1+\y)--(-0.5+\x,1+\y) node[font = \small,right = 6mm,above]{$x_6 = 1$};
			\tikzmath{\x =8.2;};
			\draw[->] (3+\x,1+\y)--(0.5+\x,1+\y) node[font = \small,right = 6mm,above]{$x_7 = 2$};
			\tikzmath{\x =0.3;};
			\draw[->] (3+\x,1+\y)--(0.5+\x,1+\y) node[font = \small,right = 6mm,above]{$x_8 = 4$};
		\end{tikzpicture}
	\end{center}
	\caption{An example of the bijection $\psi$ between $\mathcal{P}^d_n$ and $\mathcal{A}^d_n$.} \label{fig:psi}
\end{figure}

The following lemma gives another equivalent description of the map $\psi$.
\begin{lemma}\label{lem:psi-equi}
	For $d\geq 0$ and $n \geq 1$, let $P\in \mathcal{P}^d_n$ with $\omega(P) = a_1a_2\cdots a_n$ and let $x = x_1x_2\cdots x_n = \psi(P)$.
	Then for $i \in [1,n]$, we have $x_i$  is the number of active elements of $P$
	in $[a_i]$ with the convention that $[0] = \emptyset$.
\end{lemma}

\pf 
We prove the result by induction on $n$ where $n=1$ is trivial.
Assume the result for $n-1$.
By the recursive construction of $\psi$, it is easily seen that 
$x' = x_1x_2\cdots x_{n-1} = \psi(P[n-1])$.
Note that $\omega(P[n-1])=a_1a_2\cdots a_{n-1}$.
Then by the induction hypothesis,  we have $x_i$ $(1\leq i \leq n-1)$ is the number of active elements of the poset $P[n-1]$ in $[a_i]$.
Notice that $\mathrm{Act}(P)$ and $\mathrm{Act}(P[n-1])$ may differ by one element ${n-1}$.
Combining the fact that 
 $a_i\leq i-1 <n-1$ $(1\leq i \leq n-1)$,  
 we deduce that
 $P$ and $P[n-1]$ have the same active elements in $[a_i]$ $(1\leq i \leq n-1)$.
Hence $x_i$  is also the number of active elements of $P$
in $[a_i]$ for $1\leq i \leq n-1$.

It remains to prove that $x_n$ is the number of active elements of $P$
in $[a_n]$.
There are three cases for $a_n$
where the first and the third cases of (\ref{equ:psi}) are trivial.
For the case $a_n = n-1$, by the definition of $\psi$, we obtain that $x_n= \mathrm{act}(P[n-1])+1$.
Recall that 
 $n-1 \in \mathrm{Act}(P)$ in this case.
 Hence the poset $P$  contains $\mathrm{act}(P[n-1])+1 = x_n$ active 
 elements in $[n-1]$, completing the proof.
 \qed

\begin{theorem}\label{thm:psi}
	For $d \geq 0$ and $n\geq 1$, the map $\psi$ is a bijection between 
	$\mathcal{P}^d_n$ and $\mathcal{A}^d_n$.
	Furthermore,  we have $\mathrm{Act}(P) = \mathrm{dAsc}(\psi(P))$
	for all $P \in \mathcal{P}^d_n$.
\end{theorem}

\pf 
Since the sequence $\psi(P)$ encodes the construction of $P$, the map 
$\psi$ is injective.
In order to prove $\psi$ is a bijection, we need to show that 
$\psi$ is surjective, namely, $\psi(\mathcal{P}^d_n) = \mathcal{A}^d_n$.
The recursive construction of the map $\psi$ tells us that 
$x = x_1x_2\cdots x_n\in \psi(\mathcal{P}^d_n)$ if and only if 
$x' = x_1x_2\cdots x_{n-1} \in \psi(\mathcal{P}^d_{n-1})$
and $0\leq x_n \leq \mathrm{act}(\psi^{-1}(x')) +1$.
By induction on $n$ and the definition of $d$-ascent sequences, 
to prove $\psi(\mathcal{P}^d_n) = \mathcal{A}^d_n$, 
it is sufficient to show that $\mathrm{Act}(P) = \mathrm{dAsc}(\psi(P))$
for any $P \in \mathcal{P}^d_n$.

Now let us focus on the property $\mathrm{Act}(P) = \mathrm{dAsc}(\psi(P))$.
We will prove this result by induction on $n$ where $n=1$ is trivial.
Assume the result for $n-1$.
Let $P$ be a difference $d$ poset in $\mathcal{P}^d_n$ and $x = x_1x_2\cdots x_n = \psi(P)$.
Then we have $x' = x_1x_2\cdots x_{n-1} = \psi(P[n-1])$.
By the induction hypothesis, we have 
$\mathrm{Act}(P[n-1]) = \mathrm{dAsc}(x')$.
To prove $\mathrm{Act}(P) = \mathrm{dAsc}(x)$,
it is sufficient to prove $n-1\nin \mathrm{Act}(P)$ 
if and only if $n-1\nin \mathrm{dAsc}(x)$.
By the definition of  inactive elements, 
$n-1 \nin \mathrm{Act}(P)$  if and only if $a_{n}\leq a_{n-1}$ and there are at least $d$ active elements in $[a_{n}+1, a_{n-1}]$.
From Lemma \ref{lem:psi-equi}, 
this is equivalent to the fact $x_n \leq  x_{n-1} - d$, namely, 
$n-1 \nin \mathrm{dAsc}(x)$.
To conclude, we have $\mathrm{Act}(P) = \mathrm{dAsc}(\psi(P))$.
This completes the proof.
\qed

Combining Theorems \ref{thm:d-poset} and \ref{thm:psi} gives a 
new proof of Theorem \ref{thm:poset-injection}.
It should be mentioned that the map $\psi$ also applies to $d=0$.
When $d=  0$,  the map $\psi$ induces a bijection between $\mathcal{P}^0_n$ and 
ascent sequences $\mathcal{A}_n$,
where $\mathcal{P}^0_n$ can be equivalently described as the set of factorial posets $P$ on $[n]$ satisfying the following  rule: 
if there exist some $i,k$ such that $i<_P k$ and $i+1 \nless_P k$, 
then there exists some $j$ such that $j<_P i+1$ and $j \nless_P i$.
We remark that the inverse of the map $\psi$ is the map $\mathrm{po}$ in 
Theorem \ref{thm:poset-injection}
which can be verified from Lemma \ref{lem:psi-equi},
Theorem \ref{thm:psi} and  the proof of Theorem 5.4 in 
\cite{Dukes2023}.
We omit the detailed proof here.

\section{Matrices}\label{sec:matrix}
	
In this section, we give an answer to a problem posed by 
Dukes and Sagan \cite{Dukes2023} by building a direct bijection between 
a class of matrices with a certain column restriction and Fishburn matrices.

Fishburn matrices were introduced by Fishburn \cite{Fishburn3} to represent
interval orders.
The first explicit bijection $\mathrm{mx}'$ between ascent sequences and Fishburn matrices was given by Dukes and Parviainen \cite{Dukes2010}.
In order to solve a conjecture of Jel\'inek \cite{Jelinek}, 
Chen et al. \cite{Yan2019} 
constructed another bijection $\mathrm{mx}''$ between ascent sequences and Fishburn matrices.
A {\em Fishburn matrix} $A$ is an upper triangular matrix with nonnegative integers such that all rows and columns contain at least one nonzero entry.
We define the {\em weight} of a matrix $A$, denoted by $w(A)$, to be 
the sum of the entries of $A$. 
Let $\mathcal{M}_n$ denote the set of Fishburn matrices of weight $n$.
For example, we have 
\[
\mathcal{M}_3=\left\{
\begin{pmatrix} 3 \end{pmatrix},
\begin{pmatrix} 2 & 0 \\0& 1 \end{pmatrix},
\begin{pmatrix} 1 & 1 \\0& 1 \end{pmatrix},
\begin{pmatrix} 1 & 0 \\0& 2 \end{pmatrix},
\begin{pmatrix} 1 & 0&0 \\0& 1&0\\0&0&1 \end{pmatrix}
\right\}.
\]

Given a  matrix $A$, let $\mathrm{dim}(A)$ denote 
the number of rows of the matrix $A$ 
and  let $A_{i,j}$ denote the entry in the $i$-th row and 
$j$-th column of $A$.
We assume that the rows  of a matrix are numbered from top to bottom  and the columns  are numbered from left to right in which the topmost row is numbered by $1$ and the leftmost column is numbered by $1$.
A row (or column) is said to be zero if all the entries in the row (or column)
are zero.

Given a matrix $A$, we let $c_j(A)$ be the  column vector consisting of the $j$-th column of $A$.
If $c_j(A)$ is not zero, then we define $\mathrm{rmax_j}(A)$ and 
$\mathrm{rmin_j}(A)$ to be the largest and the smallest index $i$ such that $A_{i,j} > 0$, respectively.
In what follows, we always assume that 
matrices are square matrices with nonnegative integers and  contain no zero columns unless specified otherwise.

Recall that Dukes and Sagan \cite{Dukes2023} constructed a bijection $\mathrm{mx}$ between 
$d$-ascent sequences and a class of upper triangular matrices with a certain column
restriction.
When $d = 0$, $\mathrm{mx}$ induces a bijection between 
ascent sequences $\mathcal{A}_n$ and a class of matrices $\mathcal{M}_n'$ defined as follows.

\begin{definition}
	Let $\mathcal{M}'_n$ be the set of upper triangular matrices $A$  with nonnegative integers which 	satisfy the following properties:
	\begin{itemize}
		\item[$(\mathrm{Ma})$] The weight of $A$ is $n$.
	    \item[$(\mathrm{Mb})$] There exist no zero columns in $A$.
	    \item[$(\mathrm{Mc})$]  For all $1 \leq j \leq \mathrm{dim}(A)-1$,
	                           $\mathrm{rmax_{j+1}}(A) > \mathrm{rmin_{j}}(A)$.
	\end{itemize}
\end{definition}

For example, let $A$ be the matrix  shown in Figure \ref{fig:matrix}.
The  $\mathrm{rmax}$ and $\mathrm{rmin}$ values of 
each column are listed below the matrix.
It can be checked that $A$ is a matrix in  $\mathcal{M}'_{19}$.
Let $\mathcal{M}' = \bigcup_{n\geq 1}\mathcal{M}'_n$.

\begin{figure}[H]
	\begin{center}
		$
		A \; = \; \begin{pmatrix}
			2 & 1 & 3 & 2 & 1\\
			0 & 1 & 1 & 3 & 0\\
			0 & 0 & 0 & 1 & 2\\
			0 & 0 & 0 & 2 & 0\\
			0 & 0 & 0 & 0 & 0\\
		\end{pmatrix} $\\
		\vskip 10pt
		\begin{tikzpicture}[font =\small , scale = 0.4]
			\foreach \j in {1,2,3,4,5} \node at (1.35*\j,0) {$\j$};
			\foreach \j /\k in {1/1,2/1,3/1,4/1,5/1} \node at (1.35*\j,-1.5) {$\k$};
			\foreach \j /\k in {1/1,2/2,3/2,4/4,5/3} \node at (1.35*\j,-3) {$\k$};
			\tikzmath{\c = -0.4;};
			\node at (\c,0){$j$:};
			\node at (\c,-1.5){$\mathrm{rmin}_j$:};
			\node at (\c,-3){$\mathrm{rmax}_j$:};
		\end{tikzpicture}
	\end{center}
	\caption{A matrix of $\mathcal{M}'_{19}$ with its $\mathrm{rmax}$ and $\mathrm{rmin}$ values.} \label{fig:matrix}
\end{figure}

 
 \begin{problem}(\cite{Dukes2023}, Problem 8.3)\label{prob:2}
 	Find a direct bijection $\mathcal{M}'_n \rightarrow \mathcal{M}_n$ 
 	without composing $\mathrm{mx}^{-1}$ and $\mathrm{mx}'$.
 \end{problem}

We will give an answer to Problem \ref{prob:2} by 
constructing a map $\theta$ from  $\mathcal{M}'_n$ to $\mathcal{M}_n$ 
which we then show (in Theorem \ref{thm:theta}) to be a bijection. 
To this end, we need to define two transformations $\alpha$  and $\beta$ on matrices which will play essential roles in the construction of $\theta$ and its inverse $\theta'$, respectively.
For $1\leq k \leq \mathrm{dim}(A)$,  let $A[k]$ denote the submatrix of $A$
composed of entries from the first $k$ rows and first $k$ columns of $A$.

\noindent{\bf The transformation $\alpha$ } \\
Let $A$ be a matrix with $\mathrm{dim}(A) = m$ and $\mathrm{rmax}_m(A) = i$.
The matrix $\alpha(A)$ is defined as follows.

\begin{itemize}
	\item [(1)] If $i= m$, then let $\alpha(A) = A$.
	\item [(2)]	 If $i < m$,  then we construct 
	$\alpha(A)$ in the following way.       
	In $A[m-1]$, insert a new zero row between rows $i-1$ and $i$, 
	and insert a new zero column between columns $i-1$ and $i$.
	Denote by $A'$ the resulting matrix.
	Then copy the highest $i-1$ entries in the last column of $A'$ 
	to the top $i-1$ entries in the new zero column.
	Set $A''$ to be the resulting matrix.
	Finally replace the highest $i$ entries in the last column of $A''$
	with the top $i$ entries in the last column of $A$.	       
	The resulting matrix is $\alpha(A)$.
\end{itemize}

\begin{example}\label{exam:alpha}
	Consider the following two matrices:
	$$
	A=\left(
	\begin{array}{cccc}
		1 & 1 & 3 & 0 \\
		0 & 2 & 1 & 1 \\
		0 & 0 & 1 & 0 \\
		0 & 0 & 0 & 2 \\
	\end{array}
	\right);
	\quad
	B=\left(
	\begin{array}{ccccc}
		1 & 0 & 2 & 3 & 1\\
		0 & 3 & 1 & 1 & 0\\
		0 & 0 & 2 & 3 & 2\\
		0 & 0 & 0 & 2 & 0\\
		0 & 0 & 0 & 0 & 0\\
	\end{array}
	\right).
	$$
For matrix $A$, we have $\mathrm{rmax}_4(A) = \mathrm{dim}(A) = 4$, 
by using rule $(1)$ of the
transformation $\alpha$, we have   $\alpha(A) = A$.
For matrix $B$, since $\mathrm{rmax}_5(B) =3 < \mathrm{dim}(B) = 5$,
rule $(2)$ of the
transformation $\alpha$
applies and we do as follows.
Insert a zero row and zero column between rows $2$ and $3$
and columns $2$ and $3$ of $B[4]$, respectively.
We obtain the resulting matrix $B'$  shown as follows, 
with the newly inserted zeros highlighted in bold.
\[
B'=
\begin{pmatrix} 	
	1 & 0 & {\bf 0} & 2 & 3 \\
	0 & 3 & {\bf 0} & 1 & 1 \\
	{\bf 0} & {\bf 0} & {\bf 0} & {\bf 0} & {\bf 0} \\
	0 & 0 & {\bf 0} & 2 & 3 \\
	0 & 0 & {\bf 0} & 0 & 2 \end{pmatrix}.
\]
Next copy the highest $2$ entries in the last column of $B'$
to the top $2$ entries in the new zero column.
These are illustrated in red in the following matrix:
\[B''=
\begin{pmatrix} 	
	1 & 0 & {\bf \red{3}} & 2 & 3 \\
	0 & 3 & {\bf \red{1}} & 1 & 1 \\
	{\bf 0} & {\bf 0} & {\bf 0} & {\bf 0} & {\bf 0} \\
	0 & 0 & {\bf 0} & 2 & 3 \\
	0 & 0 & {\bf 0} & 0 & 2 \end{pmatrix}.
\]
Finally replace the highest $3$ entries in the last column of $B''$
with the top $3$ entries in the last column of $B$ to yield $\alpha(B)$:
\[
\alpha(B)=
\begin{pmatrix} 	
	1 & 0 & {\bf \red{3}} & 2 & {\bf \red{1}} \\
	0 & 3 & {\bf \red{1}} & 1 & {\bf \red{0}} \\
	{\bf 0} & {\bf 0} & {\bf 0} & {\bf 0} & {\bf \red{2}} \\
	0 & 0 & {\bf 0} & 2 & 3 \\
	0 & 0 & {\bf 0} & 0 & 2 \end{pmatrix}.
\]
\end{example}	

Given a  Fishburn matrix $A$ with $\mathrm{dim}(A) = m$, 
let $\mathrm{index}(A)$ denote the smallest index $i$ such that 
the $i$-th row of $A$ is zero everywhere except for entry $A_{i,m}$.
Since the only nonzero entry in the $m$-th row of $A$ is $A_{m,m}$, 
then $\mathrm{index}(A)$ is well-defined.


\begin{lemma}\label{lem:alpha}
	Let $A$ be  a matrix with $\mathrm{dim}(A) = m$ such that 
	$A[m-1]$ is a Fishburn matrix and $\mathrm{rmin}_{m-1}(A) < \mathrm{rmax}_{m}(A)$.
	Then we have that $\alpha(A)$ is a Fishburn matrix.
	Moreover we have the weight $w(\alpha(A)) = w(A)$, $\mathrm{dim}(\alpha(A)) = \mathrm{dim}(A)$, $\mathrm{rmin}_m(\alpha(A)) = \mathrm{rmin}_m(A)$
	and $\mathrm{index}(\alpha(A)) = \mathrm{rmax}_m(A)$.
\end{lemma}

\pf
We first show $\alpha(A)$ is a Fishburn matrix.
We have two cases.
To simplify notation, let $i = \mathrm{rmax}_m(A)$.
If $i = m$,  by using rule $(1)$ of 
$\alpha$, we have  $\alpha(A) = A$.
In this case, we have $A_{m,m} >0$.
Combining the fact that $A[m-1]$ is a Fishburn matrix, 
we deduce that $A$ contains no zero rows or zero columns.
It is easily seen that $A$ is an upper triangular matrix.
Hence, $\alpha(A) = A$ is a Fishburn matrix.

If $i < m$, rule $(2)$ of $\alpha$ applies.
We first claim that $\alpha(A)$ does not contain zero rows or zero columns.
Since $A[m-1]$ is a Fishburn matrix,  then by the construction of 
$\alpha$,
it is sufficient to show that  the newly added row (resp. column) is not zero.
Notice that $\mathrm{rmin}_{m-1}(A) < \mathrm{rmax}_{m}(A)$.
We have that the newly added column is not zero.
Again by the construction of $\alpha$, 
we deduce that the last entry of the  newly added row is not zero.
This proves the claim.
It is routine to check that $\alpha(A)$ is still  an upper triangular 
matrix.
This yields that $\alpha(A)$ is a Fishburn matrix.

For the second part of the lemma, we
will only prove $\mathrm{index}(\alpha(A)) = \mathrm{rmax}_m(A)$ as the other equalities can be easily
verified by the construction of the transformation $\alpha$. 
Since $A[m-1]$ is a Fishburn matrix,  we have that each row
of $A[m-1]$ contains at least one nonzero entry.
Again by the construction of $\alpha$, it can be checked that 
the first $i-1$ rows of  $\alpha(A)[m-1]$  are not zero.
Moreover, the $i$-th row of the matrix $\alpha(A)$ contains 
a unique nonzero entry $\alpha(A)_{i,m}$.
Hence we have $\mathrm{index}(\alpha(A)) = i=\mathrm{rmax}_m(A)$,
as desired.
\qed

\noindent{\bf The transformation $\beta$ } \\
Let $A$ be a Fishburn matrix with $\mathrm{dim}(A) = m$ and $\mathrm{index}(A) = i$.
The matrix $\beta(A)$ is defined as follows.

\begin{itemize}
	\item[(1)]  If $i= m$, then let $\beta(A) = A$.
	\item[(2)] 	If $i < m$,  then we construct 
	$\beta(A)$ in the following way.      
	Let $B$ be the matrix obtained from $A$ by  replacing  the highest $i$ entries in the last column of $A$
	with the top $i$ entries in the $i$-th column of $A$.	
	Delete the $i$-th row and $i$-th column of $B$, then 
	insert a zero row at the bottom and a zero column to the right of $B$.
	Let $C$ be the resulting matrix.
	Then $\beta(A)$ is the matrix obtained from $C$ by replacing  the highest $i$ entries in the last column of $C$
	with the top $i$ entries in the last column of $A$.	
\end{itemize}

For example, consider the matrix $A$ in Example \ref{exam:alpha}.
Since $\mathrm{index}(A) = \mathrm{dim}(A) = 4$, 
rule $(1)$ of the transformation $\beta$ applies and 
we obtain $\beta(A) = A$.
For  the matrix $\alpha(B)$ in Example \ref{exam:alpha}, 
since $\mathrm{index}(\alpha(B)) =3 <\mathrm{dim}(\alpha(B)) = 5$,
we obtain $\beta(\alpha(B)) = B$ by applying  rule (2) of the transformation $\beta$.


\begin{lemma}\label{lem:alpha-beta}
	Let $A$ be a matrix with $\mathrm{dim}(A) = m$ such that $A[m-1]$ is a Fishburn matrix and $\mathrm{rmin}_{m-1}(A) < \mathrm{rmax}_{m}(A)$.
	Then we have that $\beta(\alpha(A)) = A$.
\end{lemma}

\pf 
By Lemma \ref{lem:alpha}, we see that $\alpha(A)$ is a Fishburn matrix and  $\mathrm{index}(\alpha(A)) = \mathrm{rmax}_m(A)$.
Then the conclusion follows directly from the fact that  cases $(1)$ and $(2)$ in the construction of  $\beta$
correspond, respectively, to cases $(1)$ and $(2)$ of the construction of  $\alpha$.
\qed

Let $A$ be a  matrix. 
For $1 \leq k \leq \mathrm{dim}(A)$, let $\alpha_k(A)$ denote the
matrix obtained from $A$ by replacing the submatrix $A[k]$ 
with $\alpha(A[k])$ and keeping other entries in $A$ unchanged.
We are now ready to define our map $\theta: \mathcal{M}'_n\rightarrow \mathcal{M}_n$.
Given $A \in \mathcal{M}'_n$ with $\mathrm{dim}(A) = m$, define $$\theta(A) = \alpha_m\circ\alpha_{m-1} \circ \cdots \circ \alpha_1(A).$$

\begin{example}\label{exam:theta}
	Let $A$ be the matrix shown in Figure \ref{fig:matrix}.
	Then we have 
		\[	\theta(A)=\left(
	\begin{array}{ccccc}
		2 & 1 & 2 & 3 & 1\\
		0 & 0 & 3 & 1 & 0\\
		0 & 0 & 0 & 0 & 2\\
		0 & 0 & 0 & 1 & 1\\
		0 & 0 & 0 & 0 & 2\\
	\end{array} 
	\right) \in \mathcal{M}_{19},\]
	with the detailed process below.
	\begin{align*}
			A=\left(
		\begin{array}{ccccc}
			2 & 1 & 3 & 2 & 1\\
			0 & 1 & 1 & 3 & 0\\
			0 & 0 & 0 & 1 & 2\\
			0 & 0 & 0 & 2 & 0\\
			0 & 0 & 0 & 0 & 0\\
		\end{array} 
		\right)
		&\xrightarrow{\alpha_1}
		\left(
		\begin{array}{ccccc}
			2 & 1 & 3 & 2 & 1\\
			0 & 1 & 1 & 3 & 0\\
			0 & 0 & 0 & 1 & 2\\
			0 & 0 & 0 & 2 & 0\\
			0 & 0 & 0 & 0 & 0\\
		\end{array} 
		\right)	
		\xrightarrow{\alpha_2}
		\left(
		\begin{array}{ccccc}
			2 & 1 & 3 & 2 & 1\\
			0 & 1 & 1 & 3 & 0\\
			0 & 0 & 0 & 1 & 2\\
			0 & 0 & 0 & 2 & 0\\
			0 & 0 & 0 & 0 & 0\\
		\end{array} 
		\right)	 \\[5pt]	
		\xrightarrow{\alpha_3}
		\left(
		\begin{array}{ccccc}
			2 & 1 & 3 & 2 & 1\\
			0 & 0 & 1 & 3 & 0\\
			0 & 0 & 1 & 1 & 2\\
			0 & 0 & 0 & 2 & 0\\
			0 & 0 & 0 & 0 & 0\\
		\end{array} 
		\right)
		&\xrightarrow{\alpha_4}
		\left(
		\begin{array}{ccccc}
			2 & 1 & 3 & 2 & 1\\
			0 & 0 & 1 & 3 & 0\\
			0 & 0 & 1 & 1 & 2\\
			0 & 0 & 0 & 2 & 0\\
			0 & 0 & 0 & 0 & 0\\
		\end{array} 
		\right)
		\xrightarrow{\alpha_5}
		\left(
		\begin{array}{ccccc}
			2 & 1 & 2 & 3 & 1\\
			0 & 0 & 3 & 1 & 0\\
			0 & 0 & 0 & 0 & 2\\
			0 & 0 & 0 & 1 & 1\\
			0 & 0 & 0 & 0 & 2\\
		\end{array} 
		\right).
	\end{align*}
\end{example}

Let $A$ be a matrix in $\mathcal{M}'_n$ with $\mathrm{dim}(A) = m$.
In the rest of the paper, for $1\leq i \leq m$, we always denote $$A^{(i)} = \alpha_i(A^{(i-1)})$$ with the convention that $A^{(0)} = A$.
By the definition of $\theta$, we have that 
$\theta(A) = A^{(m)}$.

\begin{lemma}\label{lem:well-defined}
	Let $A$ be a matrix in $\mathcal{M}'_n$ with $\mathrm{dim}(A) = m$.
	Then for $1\leq i\leq m$, we have that $A^{(i)}$ is an  upper triangular matrix which satisfies the following properties:
	\begin{itemize}
		\item[\upshape (1)]  $w(A^{(i)}) = n$;
		\item[\upshape (2)]  $A^{(i)}[i]$ is a Fishburn matrix;
		\item[\upshape (3)]   $\mathrm{rmin}_i (A^{(i)})< \mathrm{rmax}_{i+1}(A^{(i)})$ with the convention that $\mathrm{rmax}_{m+1}(A^{(m)}) = m+1$.
	\end{itemize}
\end{lemma}

\pf 
From Lemma \ref{lem:alpha}, $\alpha$ is a weight preserving map.
Then (1) follows directly from the definition of $A^{(i)}$.
We will prove (2) and (3)  by induction on $m$, where the case $m = 1$ is trivial.
Assume the result for $m-1$.
We need to prove the result for $m$.
For (2), we need to show that $A^{(m)}[m] = A^{(m)}$ is a Fishburn matrix.
For (3), it is readily apparent that (3) holds for $i = m$.
However, the value of $\mathrm{rmax}_{m} (A^{(m-1)})$ may  change when $\text{dim}(A)$ increases from $m-1$ to $m$.
Hence, it is necessary to show that 
$\mathrm{rmin}_{m-1} (A^{(m-1)})< \mathrm{rmax}_{m}(A^{(m-1)})$.
Observe that $A[m-1]$ is a matrix in $\mathcal{M}'$.
Then by the induction hypothesis, 
$A^{(m-1)}[m-1]=A[m-1] ^{(m-1)}$ is a  Fishburn matrix.
By Lemma \ref{lem:alpha} and the construction of $\alpha_k$, 
one can verify that $\mathrm{rmin}_{m-1}(A^{(m-1)}) = \mathrm{rmin}_{m-1}(A) < \mathrm{rmax}_{m}(A) = \mathrm{rmax}_{m}(A^{(m-1)})$.
Note that $A^{(m)} = \alpha_m(A^{(m-1)}) = \alpha(A^{(m-1)})$.
Then from Lemma \ref{lem:alpha}, we derive that $A^{(m)}$ is a Fishburn matrix.
This completes the proof.
\qed

Given two positive integers $i \leq m$,
let $A$ be a matrix with $\mathrm{dim}(A) = m$ such that $A[i]$ is a Fishburn matrix.
Define $\beta_i(A)$ to be the
matrix obtained from $A$ by replacing the submatrix $A[i]$ 
with $\beta(A[i])$ and keeping other entries in $A$ unchanged.
Given $B\in \theta(\mathcal{M}'_n)$ with $\mathrm{dim}(B) = m$,
define 
$$\theta'(B) = \beta_1\circ \beta_2 \circ \cdots \circ \beta_{m}(B).$$
For example, consider the matrix $\theta(A)$  shown in Example \ref{exam:theta}, where $A$ is the matrix illustrated in Figure 
\ref{fig:matrix}.
We have $\theta'(\theta(A))  = A$.
See the procedure in Example \ref{exam:theta} backward for an illustration.

\begin{theorem}\label{thm:theta}
	For $n \geq 1$, the map $\theta$ is a bijection between 
	$\mathcal{M}'_n$ and $\mathcal{M}_n$.
\end{theorem}

\pf
Let $A$ be a matrix in $\mathcal{M}'_n$ with $\mathrm{dim}(A) = m$.
According to Lemma \ref{lem:well-defined}, we have that $\theta$ is well-defined and $\theta(A)= A^{(m)}[m] \in \mathcal{M}_n$.
By Lemmas \ref{lem:alpha-beta} and  \ref{lem:well-defined}, it is easily seen that $\theta'$ is also well-defined and 
$\theta'\circ \theta(A) = A$.
This implies that $\theta$ is injective.
Combining the fact that $\mathcal{M}'_n$ and $\mathcal{M}_n$ are equinumerous,
we conclude that  $\theta$ is a bijection between 
$\mathcal{M}'_n$ and $\mathcal{M}_n$.
\qed

From the proof of Theorem \ref{thm:theta}, 
it is easily seen that 
 the maps $\theta$ and $\theta'$ are inverses 
of each other.

\begin{remark}
	For $A \in \mathcal{M}'_n$, 
	we have showed that $\theta(A) = \mathrm{mx}' \circ \mathrm{mx}^{-1}(A)$ by induction on $\mathrm{dim}(A)$.
	Moreover, there is another direct bijection $\bar{\theta}$ between $\mathcal{M}'_n$ and 	$\mathcal{M}_n$.
	The construction of $\bar{\theta}$ is exactly the same as that of $\theta$
	except replacing $\alpha$ with $\alpha'$,
	where $\alpha'$ is defined as follows.
	
	\noindent{\bf The transformation $\alpha'$ } \\
	Let $A$ be a matrix with $\mathrm{dim}(A) = m$ and $\mathrm{rmax}_m(A) = i$.
	The matrix $\alpha'(A)$ is defined as follows.
	\begin{itemize}
		\item[(1)] If $i= m$, then let $\alpha'(A) = A$.
		\item[(2)] If $i < m$,  then we construct 
		$\alpha'(A)$ in the following way.       
		In $A[m-1]$, insert a new zero row between rows $i-1$ and $i$, 
		and insert a new zero column between columns $i-1$ and $i$.
		Denote by $A'$ the resulting matrix.
Let $T$ be the set of indices $j$ such that $j\geq i+1$ and  column $j$ of $A'$  contains  at least  one nonzero entry above row $i$. Suppose that $T=\{c_1, c_2, \ldots, c_\ell\}$ with $c_1<c_2<\cdots <c_{\ell}$.
  Let $c_0=i$. 
  For all $1\leq a\leq i$ and $1\leq b\leq \ell $, 
  copy the entry $A_{a,c_b}$ to $A_{a,c_{b-1}}$
  and replace the highest $i$  entries in the last column of $A'$ 
  with the top $i$ entries in the last column of $A$.   
  The resulting matrix is  $\alpha'(A)$.
	\end{itemize}	 
Similarly, for $A \in \mathcal{M}'_n$,  we can prove that  $\bar{\theta}(A) = \mathrm{mx}'' \circ \mathrm{mx}^{-1}(A)$ by induction on $\mathrm{dim}(A)$. 
The proofs of these results are left to the readers who are familiar with these maps.
\end{remark}

	\section*{Acknowledgments}
	 The authors are very grateful to the referees for valuable comments and suggestions
	which helped to improve the presentation of the paper.
	The work  was supported by
	the National Natural
	Science Foundation of China (11801378 and 12071440).
	

\end{document}